\documentclass[11pt]{amsart}
\usepackage{amsmath,amssymb}

\theoremstyle{plain}
\newtheorem{propn}{Proposition}[section]
\newtheorem{thm}[propn]{Theorem}
\newtheorem{lemma}[propn]{Lemma}
\newtheorem{cor}[propn]{Corollary}

\theoremstyle{definition}
\newtheorem{defn}[propn]{Definition}

\theoremstyle{remark}
\newtheorem*{rem}{Remark}

\numberwithin{equation}{section}

\begin{document}

\title{Rigidity results for Hermitian-Einstein manifolds}
\author{Stuart J. Hall}
\address{Department of Applied Computing, University of
Buckingham, Hunter Street, Buckingham, United Kingdom.}
\email{s.hall@buckingham.ac.uk}
\author{Thomas Murphy }
\address{Department of Mathematics, California State University Fullerton, 800 N State College Blvd., Fullerton, CA 92831 USA.}
\email{tmurphy@fullerton.edu  }

\thanks{ We are greatly indebted to Kouei Sekigawa, who
very kindly explained his work to us in detail, encouraged our research and gave us
constructive feedback. We also wish to thank everybody who gave us constructive comments and helped improve an earlier draft.}

\begin{abstract}
A differential operator introduced by A. Gray on
the unit sphere bundle of a K\"ahler-Einstein manifold is studied.
A lower bound for  the first eigenvalue of the Laplacian for the 
Sasaki metric on the unit sphere bundle of a K\"ahler-Einstein manifold is derived. Some rigidity theorems classifying  complex space forms amongst compact
Hermitian surfaces and the product of two projective lines amongst all K\"ahler-Einstein surfaces are then derived.
\end{abstract}

\maketitle

\section{introduction}

The main inspiration for this work is the following theorem:
\begin{thm}\label{1} (Gray \cite{gray})
Let $(M,g)$ be a closed K\"ahler manifold with constant scalar curvature.  If $M$ has non-negative
sectional curvature, then $(M,g)$ is isometric to a locally symmetric space of compact
type.
\end{thm}

 Gray also gives two proofs that $\mathbb{C}P^n$ equipped with the Fubini-Study metric is the unique K\"ahler manifold with positive bisectional curvature (a classical result of Berger \cite{berger}) as  a consequence of his method.  
  
The method Gray uses to prove this theorem is extremely interesting. He
introduces a linear differential operator $L$ on the unit sphere bundle $S(M)$
of $M$, whose coefficients are determined by the sectional curvatures of $M$. 
Thus, for example, when $M$ is positively curved $L$ turns out to be elliptic.
However $L$ is not as well-behaved if  one assumes $M$ has nonnegative bisectional curvature. 
 Gray states in \cite{gray} that he
expects the method to have further applications; we aim to make a small contribution in this direction.

Throughout the paper, we will assume all manifolds are smooth, connected, and
closed.  For any manifold $(M,g)$, let $S(M)$ denote the unit sphere
bundle of $(M,g)$, with fibre $S_p(M)$ over a point $p\in M$.
Equip $S(M)$ with the Sasaki metric $g_{sas}$. An almost complex structure will
always be denoted by $J$. When $(M,g)$ has an almost-complex structure $J$, we
will be concerned with  the holomorphic sectional curvature 
$$
H(x) = B_{xx^*}= R_{xJxxJx}.
$$
This is closely related to the  study of the unit sphere bundle because, as Berger \cite{berger}
noticed, $H$ can be viewed as a smooth function on $(S(M), g_{sas})$.

 The first non-zero eigenvalue of the scalar Laplacian is one of the most
important quantities associated to any metric.  There is a famous bound due to
Lichnerowicz for the first non-zero eigenvalue of the Laplacian of a closed
manifold with positive Ricci curvature.
It is of great interest to see if one can get similar bounds for other naturally
occurring families of metrics.  We adopt the convention that 
 the Laplacian has nonpositive eigenvalues.  Our goal is to derive a
universal lower bound for $\lambda_1$ of the Sasaki metric $g_{sas}$ on $S(M)$
using Gray's differential operator.

$(M,g)$ is \emph{normalized} when
$$
max_{x,y\in TM} \lbrace |sec(x,y)| = 1\rbrace.
$$
This is just to factor out rescaling the metric by homothety, and can always be
done if $M$ is not isometric to a flat torus.
\begin{thm}\label{thm5}
Let $(M^n,g)$ be a  closed normalized K\"ahler-Einstein manifold with $dim_{\mathbb{R}}(M) = n$  which is not isometric
to a complex space form. Then
$$
\lambda_1(S(M),g_{sas}) \geq -6(n+2).
$$
\end{thm}
Rescaling has been ruled out by the assumption $g$ is normalized. It is notable that one always has such behavior for such a naturally occurring
family of metrics.

In the special case that $(M,g)$ has positive Einstein constant (i.e. $M$ is a Fano
manifold), and in addition admits holomorphic vector fields, 
 the bound given by the above theorem is not optimal. One could pull
back an eigenfunction for the first eigenvalue of $(M,g)$ and obtain a better
lower bound on $\lambda_1(S(M),g_{sas})$. Theorem \ref{thm5} does however yield
information for the general Fano case, as well as for  the Ricci-flat and negative cases.

A second motivating result in writing this work is the following classical
result of Berger (which was instrumental in his proof of Theorem \ref{1} for the
case of positive sectional curvature):
\begin{thm}\label{bergfirst}\cite{berger}
Let $(M,g)$ be a K\"ahler manifold of dimension real dimension $n=2N$. Then $$\int_{S_p(M)} H  \omega_2 =
\frac{s}{(N)(N+1)}Vol(\mathbb{S}^{n-1})$$ for all $p\in M$.
\end{thm}
The term $\omega_{2}$ in the above result comes from the splitting of the of 
canonical volume form of the Sasaki metric on $S(M)$ $\omega$ as
$\omega=\omega_{1}\wedge \omega_{2}$, where $\omega_{1}$ is the volume form
for$(M,g)$ and $\omega_{2}$ is the standard volume form for the sphere with
radius 1. Finally $s$ denotes the scalar curvature of $g$.
The classification of K\"ahler-Einstein manifolds with non-zero holomorphic
sectional curvature is a notable open problem; for example $\mathbb{C}P^N\times
\mathbb{C}P^N$ with the standard metric has positive holomorphic sectional curvature.  To
our knowledge this theorem of Berger gives us the best-known rigidity result in
this direction, namely
it follows immediately from his result that a  K\"ahler-Einstein manifold has $H
\leq   \frac{s}{(N)(N+1)}$  (or $\geq$) if, and only if $H = \frac{s}{(N)(N+1)}$
and consequently $(M,g)$ is a complex space form.

We can extend this to any almost-Hermitian manifold as follows:
\begin{thm}\label{bergergeneralization}
Let $(M,g,J)$ be an almost-Hermitian manifold of real dimension $n=2N$. Then
$$
\int_{S_p(M)} H  \omega_2 = \frac{3s^*  + s}{4(N)(N+1)}Vol(\mathbb{S}^{n-1})
$$
for all $p\in M$.
\end{thm}
We define $R^*(x,y)$, the $*$-Ricci curvature, in the following way. Then
$$
R^*_{ij} : = \sum_{a=1}^n R_{a{i}^*j{a}^*}
$$
is the star Ricci tensor, and  its trace $\sum_i R^*_{ii} = \sum_{i,a}
R_{ai^*ia*}$ is  the $*$-scalar curvature $s^*$. We caution the reader
that this definition  is slightly different to the usual one in the literature:
our definition agrees with the usual definition in the K\"ahler case.  We use this convention as it is more convenient to express our
results and it gives the correct generalization of Berger's result.

We have to use  different ideas to Berger, who heavily relies on the
K\"ahler identities and local calculations in the curvature tensor to prove his
result.
Of course in the K\"ahler case $s^*= s$ and so we recover his result. The
technique we use to calculate this identity arose from Sekigwawa and Sato's work
 \cite{ss} extending  the study of
Gray's differential operator to the nearly K\"ahler case. The main application of this estimate is that, combined with the work
of Apostolov, Davidov and Muskarov \cite{apostolov}  we obtain the following
rigidity result for closed Hermitian surfaces;
\begin{cor}
Let $(M^4,g,J)$ be a closed Hermitian surface Then $$H\leq   \frac{3s^* +
s}{24}$$ (or $\geq$) if, and only if, we have equality and 
$(M,g,J)$ is isometric to $\mathbb{C}P^2$, $\mathbb{C}^2/\Gamma$ or
$\mathbb{C}H^2/\Gamma$ equipped with their standard symmetric space metrics.
\end{cor}

Finally we  can use some of Gray's ideas to  prove the following result characterizing $\mathbb{C}P^1\times\mathbb{C}P^1$ amongst the Fano K\"ahler-Einstein surfaces. We view $H_{av}$ and $H^{max}$ as functions on a K\"ahler-Einstein surface $M$, where $H^{max}(p)$ is defined to be the maximum holomorphic sectional curvature at $p\in M$. $H_{av}$ denotes the average holomorphic sectional curvature, which by Theorem \ref{bergfirst} is constant. 

\begin{thm} \label{4} Let $(M^2, g)$ be a closed K\"ahler-Einstein surface with positive scalar curvature. Then $H_{av} = \frac{2}{3}H_{max}$ if, and only if, $M$ is isometric to $\mathbb{C}P^1\times \mathbb{C}P^1$ with the product metric. 
\end{thm}

\section{Gray's differential operator}

In this section we review the techniques Gray used to prove Theorem \ref{1},
remarking  that his entire construction generalizes to the case of $J$ being an almost-Hermitian structure.  We follow the convention that $X,Y,Z \in \Gamma(TM)$ are smooth vector
fields, and $x,y,z$ denote tangent vectors at $T_pM$. We will be brief and not fully explain all the theory behind the $L$ operator: further details are available in \cite{gray}, \cite{ss}, as well as some short notes available on the second author's website.

We define the Riemannian curvature tensor  as
$$
R_{WXYZ} = g(\nabla_{[W,X]}Y - [\nabla_W, \nabla_X]Y,Z),
$$
and  set $$sec(x,y) = \frac{R(x,y,x,y)}{\|x\|^2\|y\|^2 - g(x, y)^2}$$ to be the
sectional curvature of the plane spanned by $x,y\in T_pM$ . We will occasionally
write $R(W,X,Y,Z)$ for $R_{WXYZ}$ to make calculations easier to read. 
For $x\in S_p(M)$ take an orthonormal bases $\lbrace e_1, \dots, e_n\rbrace$ of
$T_pM$ with the convention that $x=e_1$.  Then $S_p$ is the fibre of the sphere
bundle ${S(M)\rightarrow M}$ over $p\in M$. Equip $S(M)$ with the Sasaki metric,
and for
$X\in \Gamma(TM)$ denote by $X^h$ (resp. $X^v$) the horizontal (resp. vertical)
lift. Then $\lbrace e_i^h, e_i^v \rbrace$ form an orthonormal basis of
$T_{(p,x)}S(M)$. Denote by $(y_2, ....y_n)$ the corresponding system of normal
coordinates defined on a neighbourhood of $x$ in the sphere $S_p$, and let
$(x_1, \dots, x_n)$ denote the normal coordinates corresponding to $(e_i^h)$.
Set $y_{\alpha}(x) = 0$.

For each $y\in T_pM$ the tangent space $T_y(T_p(M))$ is identified with $T_{p}M$
by means of parallel translation. Under this identification we write
$\frac{\partial}{\partial u_{i}}$ to correspond to $e_{i}.$
Set $h_{ij} = R_{i x j x} = R_{i1j1}.$

\begin{defn} Gray's $L$ operator is defined in  the local normal coordinates
$\{x_{i},y_{\alpha}\}$ by
$$
L_{(p,x)} := \bigg\lbrace \sum_{i=1}^n \frac{\partial^2 }{\partial x_{i}^2} +
\frac{1}{2} \sum_{i ,j = 2}^n
h_{ij} \frac{\partial^2}{\partial y_{i}\partial y_{j}}\bigg\rbrace_{(p,x)}
$$
\end{defn}

Viewing $H(x)$ as a function on $S(M)$, Gray proves that $L(H) = 0$ when $(M,g)$ is K\"ahler-Einstein.

\begin{lemma}[Lemma 5.3 and 7.1 in \cite{gray}] \label{33}
Let $(M,g, J)$ be an almost-Hermitian Einstein manifold. Then the following holds:
\begin{enumerate}
\item $L(H^2) = 2\| grad^h H\|_{(p,x)}^{2} + R_{x, \eta(grad^vH(x)),x, \eta(grad^vH(x))}$.\\
\item $\int_{S(M)} L(H^2) \omega = 0$.
\end{enumerate}
\end{lemma}
\begin{rem}
Note that (1) clarifies two points of Lemma 5.3 of \cite{gray}. The first is
a slight abuse of notation: the vector $grad^{v}H$ is identified with the vector
in $T_{p}M$ that lifts to $grad^{v}H$.   More importantly, the curvature term
should actually be $$R_{x, \eta(grad^vH(x)),x,\eta(grad^vH(x))}.$$ 
\end{rem}

For later, the following result will also be required (see \cite{ss});
\begin{lemma}
For $x\in S(M)$,
\begin{align} 
grad^h(H)(x) &= \sum_{i=1}^n \langle (\nabla_{e_i}R)(x,Jx),x,Jx\rangle e_i^h + 2\langle R(x,(\nabla_{e_i}J)x)x,Jx \rangle\\
grad^v(H)(x) &= 4\sum_{i=2}^n \langle R(x,Jx)x, Je_i\rangle e_i^v.
\end{align}
\end{lemma}

\section{Proof of Theorems \ref{thm5} and \ref{bergergeneralization}}

The following fact is well-known.

\begin{propn} \label{lap}
Let $\mathbb{R}^n$ denote Euclidean space and $f$ a homogeneous polynomial of
degree $r \geq 1$ on $\mathbb{R}^n$. Then
$$
\int_{\mathbb{S}^{n-1}} (\mathbb{D} f) \omega_2 =
r(n+r-2)\int_{\mathbb{S}^{n-1}} f|_{\mathbb{S}^{n-1}} \omega_2.
$$
where $\mathbb{D}$ is the Laplace-Beltrami operator of $\mathbb{R}^n$ and
$\omega_2$ denotes the volume element of the round sphere $\mathbb{S}^{n-1}$
with sectional curvature $1$.
\end{propn}

Now we give the proof of Theorem \ref{thm5}.
\proof Let $\Lambda$ denote the Einstein constant of $(M,g)$.  Let $G(x) = grad^{v}H(x)$. Now Lemma \ref{33}
and the pointwise estimate on the curvature norm together imply
\begin{align*}
 2\int_{S(M)}\|grad^hH\|^2_{(m,x)} \omega =& -\int_{S(M)} sec(x, G(x))\|G(x)\|^2
\omega\\
 \leq&  \int_{S(M)} \|G(x)\|^2\omega \\
\end{align*}
where in the last line we use the fact $(M,g)$ is normalized.

This gives the bound
$$
\int_{S(M)}\|grad(H)\|^2\omega \leq \frac{3}{2}\int_{S(M)} \|
grad^v(H)(x)\|^2\omega.
$$
The idea is to use the Raleigh quotient with $H$ as a test function to estimate
$\lambda_1(S(M), g_{sas})$. $H$ is never constant by assumption, so we can
always normalize and use $H$ as a test function on the sphere bundle.  It
remains to estimate
$$
\int_{S(M)} (H -H_{av})^2 \omega= \int_{S(M)} (H^2) \omega -
Vol(S(M))\int_{S(M)} (H_{av}^2) \omega.
$$

Let us define the functions $F,f$ on $T_pM$, for $p$ fixed, by setting
$$
F(v)= R(v,Jv,v,Jv), f = F^2
$$
for $v\in T_pM$. Writing $v = \sum_{i}v_ie_i$,
\begin{equation}\label{hform}
F(v) = \sum_{i, j,k,l\geq 1} R_{ij^*jl^*} v_iv_jv_kv_l.
\end{equation}
By definition of $F$ and $f$, $F|_{S_p} = H$ and $f|_{S_p} = H^2$. This gives
$$
grad(F) = 4\sum_{i, j,k,l\geq 1} R_{i{j}^*jl^*}v_iv_jv_lR_{ij^*kl^*}e_k^v.
$$

Hext we compute $$\mathbb{D} H(v) = 4\sum_{ij}(3R_{ij} + R_{ji})v_iv_j =
16\Lambda \|v\|^2$$ for $v\in T_pM$. In particular, for $(p,x)\in S(M)$ we
obtain
$$
\mathbb{D} H(p,x) = 16\Lambda.
$$
In a similar fashion

\begin{align*}
grad(F)_{(p,x)} &= 4\sum_{l > 1}R_{xJxxJe_l} e_l^v + 4H(x)e_i^v\\
&= grad^v(H)(x)_{(p,x)} + 4H(x)e_1^v.
\end{align*}

Applying Propostion \ref{lap} to $f$ yields that
$$
\int_{\mathbb{S}^{n-1}} (\mathbb{D}f)\omega_2 =
8(n+6)\int_{\mathbb{S}^{n-1}}H^2\omega_2.
$$

But we have that
\begin{align*}
\mathbb{D} f|_{S_p(M)} =& 2\bigg( \|grad^vH\|^2 + 16H^2 + H\mathbb{D} H\bigg)\\
&= 2\bigg(\|grad^vH\|^2 + 16H^2 + 16\Lambda H\bigg)\\
\end{align*}

Rearranging, this yields
\begin{align*}
\int_{S_p(M)}(H - H_{av})^2\omega_2 = \frac{1}{4(n+2)} \int_{S_p(M)}&
\|grad^v(H)\|^2 \\
&+ \bigg( \frac{4}{(n+2)}\Lambda- H_{av}\bigg)H_{av}Vol(S_p)\omega_2.
\end{align*}

But, again using Proposition  \ref{lap}, $H_{av} = \frac{4}{n+2}\Lambda$, so
$$
\int_{S_p(M)}(H - H_{av})^2\omega_2 = \frac{1}{4(n+2)}\int_{S_p(M)}\|
grad^v(H)\|^2\omega_2.
$$
The result follows. \endproof

\begin{rem}
In many known Ricci-flat examples, Hans-Joachim Hein has pointed out to us
that on a noncompact component of the moduli space of Ricci flat metrics can admit deformations through Ricci-flat metrics which allow one to get a  lower estimate for the first non-zero eigenvalue of 
the Laplacian of $(M,g)$, and thus pulling the corresponding eigenfunctions back to $S(M)$  a  better estimate than our work. 
Of course, such deformations do not generally exist. Our bound, in contrast,  holds uniformly with a precise constant over the whole moduli space.
\end{rem}

Next is the proof of Theorem \ref{bergergeneralization}.

\proof This follows from a similar calculation to the proof of Theorem
\ref{thm5}:  take again the function $F$ defined on $T_pM$ for some $p\in M$.
Then from Equation (\ref{hform})
$$\mathbb{D} H(v) = 4\sum_{a,i,j = 1}^n \bigg( R_{aa^*ij^*} + R_{ai^*aj^*} +   R_{ai^*ja^*}\bigg) v_iv_j$$

But $$\int_{\mathbb{S}^{n-1}} (v_iv_j) \omega_2 = 0$$ if $i\neq j$. Similarly,
it is easy to see that
$$\int_{\mathbb{S}^{n-1}} (v_i^2) \omega_2 = \frac{1}{n}Vol(\mathbb{S}^{n-1})$$
for all $i$. Then from Proposition 5.1, 
\begin{align*}
\frac{(n)(n+2)}{Vol(\mathbb{S}^{n-1})} \int_{\mathbb{S}^{n-1}} H \omega_2 =& \sum_{a,i,j = 1}^n \bigg( R_{aa^*ij^*} + R_{ai^*aj^*} +   R_{ai^*ja^*}\bigg)\bigg|_{i=j}\\
=& \sum_{i,j = 1}^n \bigg(R^*_{i^*j^*} + R^*_{ji} + R_a^*{i^*j^*} + R^*_{ij}\bigg)\bigg|_{i=j}\\
=& \sum_{i =1}^n R^*_{i^*i^*} + 2R^*_{ii} + R_{i^*i^*}.\\
\end{align*}
Since $\sum_{i=1}^n R_{i^*i^*} = \sum_{i=1}^n R_{ii} = s$ and $\sum_{i =1}^n R^*_{i^*i^*} =\sum_{i =1}^n R^*_{ii} = s^*$, the result follows.
\endproof

The corollary then follows immediately from work of Apostolev et. al \cite{apostolov}.

\section{Proof of Theorem \ref{4}}

\begin{proof} From the assumption $H_{av} = \frac{2}{3}H^{max}$ at every point $p$, we see that $H$ achieves its maximum at every point $p\in M$. At an arbitrary point $p$, choose $e_1\in T_pM$ so that $H(e_1) = H^{max}$ and then choose a normal frame as above.  Tracing the Einstein equation we obtain
$$
H_1 + B_{12} = \Lambda,
$$
where $B_{12} = R_{11^*22^*}$ is the bisectional curvature. From the assumption, $B_{12}=0$. The function $f = R_{11^*11^*}$ is constant on $M$, so at $p$   
\begin{align*}
0  = \Delta f(p) =& \sum_{i = 1}^{N} \nabla^2_{ii} f =  \sum_{i=1}^N \nabla^2_{ii} R_{11^*11^*}(p)\\
 = &\sum_{i=1}^N (\nabla^2_{ii} R)_{11^*11^*}(p)  =  \Delta^{\mathfrak{h}}(H) (p,x) .
\end{align*}
Thus, via Equation (5.3) in \cite{gray}
$$
0 \geq \Delta^{\mathfrak{h}}(H) (p,x) = (H_1 - B_{12})B_{12} - 4R_{1212}R_{12^*12^*} + 4R_{1212^*}^2.
$$
Here the fact that $R_{1211^*}$ and $R_{22^*21}$ vanish at $p$ is freely used. This follows from the fact that sectional curvatures of orthogonal planes are equal for an Einstein four manifold and the formula for $grad^vH(p,x)$.  Since  $B_{12}=0$,  $R_{1212} $ and $R_{12^*12^*}$ have the same sign if they are nonzero. But $$B_{12} = R_{1212}  +R_{12^*12^*}.$$ Hence they are both zero, and so therefore is also $R_{1212^*}$.  Repeating the above argument at every point, we see that on $(M,g)$ we have
$$
R_{1212} = R_{12^*12^*} = R_{1211^*} = R_{22^*21} = 0.
$$ 
 Let us denote by $f_i$ the usual basis of $T_0(\mathbb{C}P^1\times \mathbb{C}P^1)$ equipped with the inner product induced from the standard metric (so $f_1$ and $f_{1^*}$ span the first factor)  and  consider the map $\Phi: e_i \rightarrow f_i$ identifying the tangent space $T_qM$ with $T_0(\mathbb{C}P^1\times \mathbb{C}P^1)$.

Now pick any $y_1, y_2\in T_qM$. Then via polarization one may express $sec(y_1,y_2)$ in terms of holomorphic sectional curvatures. But the holomorphic sectional curvature of a tangent vector $\zeta$, expressed in terms of the basis $e_i$, are given by the same calculation as calculating $H(\Phi(\zeta))$ in the $f_i$ basis.
 Setting $\zeta = \sum a_ie_i = \sum a_if_i$, this can be seen from the calculation
\begin{align*}
H({\zeta}) =& R_{ \sum a_ie_i, J( \sum a_ie_i),  \sum a_ie_i,  J(\sum a_ie_i)}\\
=& \sum_{i,j,k,l} a_ia_ja_ka_l R_{e_ie_{j^*}e_ke_{l^*}}\\
=& (a_1^4 + a_{1^*}^4)H_1 + (a_2^4 + a_{2^*}^4)H_2\\
=& H(\Phi(\zeta)).
\end{align*}
 Therefore $sec(y_i, y_j)$is non-negative, because $\mathbb{C}P^1\times \mathbb{C}P^1$ with its standard metric has nonnegative sectional curvature. This is independent of the point $p\in M$. Hence $M$ has non-negative sectional curvatures, and then Theorem \ref{1} implies the result.  \end{proof}

\end{document}